\documentclass[12pt]{amsart}
\usepackage[utf8]{inputenc}
\usepackage{amssymb}
\usepackage{hyperref}
\usepackage[final]{showkeys}

\input xy
\xyoption{all}

\theoremstyle{definition}
\newtheorem{mydef}{Definition}[section]

\newtheorem{theorem}[mydef]{Theorem}
\newtheorem{cor}[mydef]{Corollary}
\newtheorem{claim}[mydef]{Claim}

\newcommand{\ba}{\bar{a}}
\newcommand{\bb}{\bar{b}}

\newcommand{\cf}[1]{\text{cf} (#1)}
\newcommand{\seq}[1]{\langle #1 \rangle}
\newcommand{\rest}{\upharpoonright}

\newbox\noforkbox \newdimen\forklinewidth
\forklinewidth=0.3pt \setbox0\hbox{$\textstyle\smile$}
\setbox1\hbox to \wd0{\hfil\vrule width \forklinewidth depth-2pt
 height 10pt \hfil}
\wd1=0 cm \setbox\noforkbox\hbox{\lower 2pt\box1\lower
2pt\box0\relax}

\title{No maximal models from looking down}
\date{\today}

\author{Will Boney}
\email{wboney@math.harvard.edu}
\address{Mathematics Department\\Harvard University\\Cambridge, MA, USA}
\date{\today\\
This material is based upon work done while the author was supported by the National Science Foundation under Grant No. DMS-1402191.} 

\begin{document}

\maketitle

\section{Introduction}

In \cite{sh893}, Shelah proved the following dichotomy for AECs with arbitrarily large models.

\begin{theorem}[\cite{sh893}.3.16]
Let $K$ be an AEC with arbitrarily large models.  Then there is a club of cardinals $\bold{C}$ such that at least one of the following holds:
\begin{itemize}
	\item for every $\lambda \in \bold{C}$ of cofinality $\omega$, there are at least $\lambda$-many nonisomorphic models of size $\lambda$.
	\item for every $\lambda \in \bold{C}$ of cofinality $\omega$, every $M \in K_\lambda$ lives inside an EM model.
\end{itemize}
\end{theorem}

The second conclusion means that there is some blueprint $\Phi$ that is proper for linear orders (see the discussion in Section \ref{outline-sec}) such that $M \prec EM_\tau(I, \Phi)$ for every linear order $I$.  In particular, this means that $M$ has extensions to models of all sizes, which is how Shelah phrases the conclusion; our statement is stronger and is proved by Shelah.

This result might suggest a new dividing line in the classification theory for AECs.  However, there is no argument given to suggest that these properties contradict each other; indeed, there are many elementary classes that satisfy both of these.  In this note, we further the argument against this being a dividing line by showing that a stronger version of the second result holds in any AEC with amalgamation; this is Theorem \ref{main-thm}.

Section \ref{outline-sec} gives the necessary background on blueprints and EM models and outlines Shelah's argument from \cite{sh893}.  Section \ref{amalg-sec} gives the proof that amalgamation can be used to replace the ``not too many models'' hypothesis.  We assume that the reader has a general familiarity with AECs; the standard references are Baldwin \cite{baldwinbook}, Grossberg \cite{ramibook}, and Shelah \cite{shelahaecbook}.

\section{Outline of Shelah's Argument} \label{outline-sec}

A key tool in the study of AECs with arbitrarily large models are EM\footnote{Stands for Ehrenfeuct-Mostowski} models and blueprints that are proper for linear orders.  

A \emph{blueprint} or a template proper for linear orders, usually denoted $\Phi$ or $\Psi$, is formally a set of complete quantifier-free types, one of each arity $n < \omega$, in some language $\tau_\Phi$ (or $\tau_\Psi$), that cohere in the following way: if $m < n$ and $\seq{n_i : i < m}$ is an increasing subsequence of $\seq{i : i < n}$, then
$$x_0, \dots x_{n-1} \vDash p_{n} \implies x_{n_0}, \dots, x_{n_{m-1}} \vDash p_m$$
and that $x_0 \neq x_1$. The name blueprint is used because these types give all the information needed to construct a $\tau_\Phi$ structure: for any linear order $(I, <_I)$, there is a unique up to canonical isomorphism $\tau_\Phi$-structure that is generated by $I$ such that every increasing $n$-sequence from $I$ realizes the type $p_n$.  We call this structure $EM(I, \Phi)$.  This construction is familiar to first-order model theorists.

Given an AEC $(K, \prec_K)$ in a language $\tau$ and a blueprint $\Phi$ such that $\tau \subset \tau_\Phi$, one can ask the following two questions:
\begin{itemize}
	\item given any linear order $I$, $EM_\tau(I, \Phi) \in K$?; and
	\item given any linear orders $I \subset J$, $EM_\tau(I, \Phi) \prec_K EM_\tau(J, \Phi)$?
\end{itemize}
where $EM_\tau(I, \Phi) := EM(I, \Phi) \rest \tau$.  We want to restrict our attention to such blueprints for which the answers to these questions are yes, so we let $\Upsilon^{or}[K]$ (often we will drop the ${}^{or}$) denote the class of all such blueprints.  At this point, the reader doesn't need to worry about the actual construction of $EM_\tau(I, \Phi)$ from $I$ and can instead think of blueprints in $\Upsilon^{or}[K]$ as a map from linear orders to elements of $K$ with the above properties\footnote{Although it's not clear that every such map actually comes from a blueprint}.

The question of whether $\Upsilon[K]$ is nonempty is equivalent to whether $K$ has arbitrarily large models, which in turn is equivalent to whether $K$ has models of sizes larger than every cardinal less than $\beth_{(2^{LS(K)})^+}$.  This analysis comes from Shelah's Presentation Theorem and Morley's Omitting Types Theorem (see \cite[Theorem 4.15]{baldwinbook} and \cite[Appendix A]{baldwinbook}, respectively).  We have so far not mentioned the \emph{size} of the blueprints, but the proof of this fact shows that there is a member of $\Upsilon[K]$ of size $LS(K)$.  If we need to emphasize the size of the blueprint, we will add a subscript, writing $\Upsilon_\kappa[K]$ or $\Upsilon^{or}_\kappa[K]$.

Given $M \in K$, we define the AEC of strong extensions of $M$ called $K_M$.  Formally, we expand the language to have a constant for each member of $M$ and require that this map is an isomorphism from $M$ to a strong substructure of the model.  Then the question of whether $K$ has no maximal models is equivalent to $\Upsilon[K_M]$ being nonempty for every $M \in K$.

We are now ready to outline Shelah's argument from \cite{sh893}.  Suppose that $M$ is of large size of countable cofinality such that
$$\lambda < \|M\| \implies \beth_{(2^\lambda)^+} < \|M\|$$
Thus, we can write $M$ as the $\prec_K$-increasing union of $\seq{M_n : n < \omega}$ such that $M_{n+1}$ is large enough to witness that $K_{M_n}$ has arbitrarily large models; write $\lambda_n:= \|M_n\|$.  Thus there is $\Phi_n \in \Upsilon[K_{M_n}]$.  If we could arrange that the $\Phi_n$ were increasing in a ``nice'' way, then this would be enough: we could define the union of them and this would be in $\Upsilon[K_M]$.  However, a priori, it is not always clear how to do this\footnote{In the next section, we use amalgamation to get a good enough approximation to do this.}.

In order to arrange this, Shelah introduces a weak notion of embedding of $EM_\tau(\lambda, \Phi)$ into $N$ over $M$ when $\Phi \in \Upsilon[K_M]$.  This is given two equivalent formulations, \cite[Definition 3.1]{sh893} as a winning strategy in a game and \cite[Definition 3.2]{sh893} as a tree of a certain depth, but we omit the details here.  We will write 
$$EM_\tau(\lambda, \Phi) \leadsto_M^{\cf \gamma} N$$
to mean that there is a direct witness to for $(N, M, \lambda, \|M\|, \gamma, \Phi)$ (in Shelah's terms) to emphasize that this is a weak notion of embedding.  This is not quite enough either, so Shelah introduces the notion of an \emph{indirect} witness, which we will write $EM_\tau(\lambda, \Phi) \leadsto_M^{\cf \gamma} N$ and means that there is an extension $\Psi$ of $\Phi$ (in the sense of \cite[Definition 2.12.(2A)]{sh893}) such that $EM_\tau(\lambda, \Psi) \leadsto_M^{\cf \gamma} N$

To construct $\Phi_n$\footnote{This construction is essentially the proof of \cite[Claim 3.12]{sh893}}, \cite[Claim 3.6.(1)]{sh893} establishes the base case.  For the inductive step, we start with $EM_\tau(\lambda_{n+4}, \Phi_n) \leadsto_{M_n}^{\cf \gamma} M$.  The structure assumption--$I(\lambda, K) \leq \lambda$--is used in \cite[Claim 3.14]{sh893} to increase the size of the $EM$ model, i. e. to get $EM_\tau(\lambda_{n+5}, \Phi_n) \leadsto_{M_n}^{\cf \gamma} M$.  This increase in size proves the existence, via \cite[Claim 3.11]{sh893}, of a $\Phi_{n+1}$ that is an extension of $\Phi_n$ in the sense of \cite[Definition 2.9.(1)]{sh893} such that $EM_\tau(\lambda_{n+5}, \Phi_{n+1}) \leadsto_{M_n}^{\cf \gamma} M$.  This allows the induction to continue.

At the end of the construction, the $\Phi_n$'s are increasing in a nice enough way that they have a natural union $\Phi$ and this union is in $\Upsilon[K]$.  Moreover, since $\Phi_n$ always contains a canonical copy of $M_n$, the union $\Phi$ always contains a canonical of $\cup_n M_n = M$.  Thus, $\Phi \in \Upsilon[K_M]$, as desired.

The reader is likely wondering why Shelah's argument uses $\lambda_{n+4}$ and $\lambda_{n+5}$ in place of $\lambda_{n+1}$ and $\lambda_{n+2}$.  The author believes this replacement to hold (especially examining the relevant hypotheses of \cite[Claim 3.9, Claim 3.11, and Theorem 3.14]{sh893}), but decided to leave Shelah's original choices intact to avoid the possibility of introducing errors.

One of the amazing things about Shelah's result is that it provides  ``downward looking'' criteria for $\Upsilon[K_M]$ to be non-empty.  If one understands the entire AEC $K$, it is easy to determine if $\Upsilon[K_M]$ is nonempty; one looks to see if $M$ has arbitrarily large extensions.  Shelah's methods, however, allow one to determine if $\Upsilon[K_M]$ is nonempty merely by looking at $K_{<\|M\|}$, the AEC below $M$.  This is very useful in the many inductive arguments that come up in the study of AECs.

\section{Amalgamation Instead of Few Models} \label{amalg-sec}

Now we prove the main result.  For ease, we define the familiar ``Hanf number function'' $h_\alpha(\kappa)$ by induction so
\begin{itemize}
	\item $h_0(\kappa) = \kappa$;
	\item $h_{\alpha + 1}(\kappa) = \beth_{(2^{h_\alpha(\kappa)})^+}$; and
	\item $h_\delta(\kappa) = \sup_{\alpha < \delta} h_\alpha(\kappa)$ for $\delta$ limit.
\end{itemize}

\begin{theorem} \label{main-thm}
Let $K$ be an AEC such that $K_{<\lambda}$ has amalgamation and $\lambda = h_\delta(LS(K)$ for $\delta$ limit. If $M \in K_\lambda$, then $\Upsilon[K_M]$ is nonempty.
\end{theorem}

We have dropped the explicit assumption that $K$ has arbitrarily large models, but it is already implied by the fact that $K$ has models of size $h_\delta(LS(K))$; it is already implied by $K$ having models of size $h_1(LS(K))$.  Also, full amalgamation in $K_{<\lambda}$ is not necessary.  As will be seen in the proof, it is enough to find a resolution of $M$ of the appropriate sizes where each $M_i$ is an amalgamation base.

{\bf Proof:}  Since $M$ is of size $h_\delta(LS(K))$, we can write $M$ as the increasing union of $\seq{M_i \in K : i < \delta}$ such that $\|M_i\| = h_i(LS(K))$.  Thus, $M_{i+1}$ witnesses that $\Upsilon[K_{M_i}]$ is nonempty.

First, we are going to construct $\Psi_i \in \Upsilon_{h_i(LS(K))}[K_{M_i}]$ such that $EM_\tau(\omega, \Psi_i)$ is embeddable into $EM_\tau(\omega, \Psi_{i+1})$.  The size of the blueprints is important because we only have amalgamation on $K_{<\lambda}$ and we need to ensure that the size of the $EM$ models is small.  Let $\Psi_0 \in \Upsilon_{h_0(LS(K))}[K_{M_0}]$ be arbitrary.  If we have $\Psi_i$, first amalgamate $M_{i+1}$ and $EM_\tau(\omega, \Psi_i)$ over $M_i$ as follows
\[
\xymatrix{
M_{i+1} \ar[r] & M_{i+1}^*\\
M_i\ar[r] \ar[u] & EM_\tau(\omega, \Psi_i) \ar[u]_{f_i}
}
\]
By amalgamating $M_{i+1}^*$ and $M_{i+2}$ over $M_{i+1}$, we know that $M_{i+1}^*$ has arbitrarily large extensions, so there is $\Psi_{i+1} \in \Upsilon_{h_{i+1}(LS(K))}[K_{M_{i+1}^*}]$; note that $\Upsilon[K_{M_{i+1}^*}] \subset \Upsilon[K_{M_{i+1}}]$.  Thus, $f_i: EM_\tau(\omega, \Psi_i) \to M_{i+1}^* \prec_K EM_\tau(\omega, \Psi_{i+1})$ as desired.  

Now suppose that $i$ is limit and $\seq{\Psi_j : j < i}$ is defined.  Then set $EM_{<i}^*$ to be the direct limit of the $EM_\tau(\omega, \Psi_i)$ with the embeddings.  Then $M_i \prec EM_{<i}^*$ and, again by amalgamating with $M_{i+1}$, we know that $\Upsilon[K_{EM_{<i}^*}]$ is nonempty.  Thus, any $\Psi_i \in \Upsilon[K_{EM_{<i}^*}]$ is as desired.

Second, we show that we can extend the $\Psi_i$'s so that the above is true for any linear order instead of just $\omega$.  Set $\Psi_0' = \Psi_0$ and for $i > 0$, set $\Psi_i'$ such that, for any linear order $I$,
$$EM_\tau(I, \Psi') = EM_\tau(\omega^i \times I, \Psi)$$
\cite[Definition 2.9.(4)]{sh893} would write this as (a particular case of) $\Psi_i \leq^{\otimes} \Psi_i'$.

In order to show that the $\Psi_i'$'s work, it suffices to prove the following claim.

\begin{claim}
If $EM_\tau(\omega, \Phi) \prec_K EM_\tau(\omega, \Psi)$, then for any linear order $I$, $EM_\tau(I, \Phi)$ is $K$-embeddable in $EM_\tau(\omega \times I, \Psi)$.
\end{claim}
From the hypothesis, we know that $EM_\tau(1, \Phi) \prec_K EM_\tau(\omega, \Psi)$, where $1$ is the one element linear order consisting of $0$.  This means there is a term $\sigma \in \tau_\Psi$ and $n_0 < \dots < n_{k-1} < \omega$ such that $0 = \sigma(n_0, \dots, n_{k-1})$.

Since each linear order is the direct limit of its finite subsets, we can write the following.
\begin{eqnarray*}
EM_\tau(I, \Phi) &=& \varinjlim_{\ba \subset \bb \in [I]^{<\omega}} \left(EM_\tau(\bb, \Phi), f_{\ba, \bb}\right) \\
EM_\tau(\omega\times I, \Psi) &=& \varinjlim_{\ba \subset \bb \in [I]^{<\omega}} \left(EM_\tau(\omega \times \bb, \Psi), g_{\ba, \bb}\right)
\end{eqnarray*}
where $f_{\ba, \bb}: EM_\tau(\ba, \Phi) \to EM_\tau(\bb, \Phi)$ is generated by the inclusion map from $\ba$ to $\bb$ and $g_{\ba, \bb}$ is similar from the inclusion map taking $\omega\times \ba$ to $\omega \times \bb$.  Then, for any $\bb \in [I]^{<\omega}$, we can define $h_{\bb}:EM_\tau(\bb, \Phi) \to EM_\tau(\omega \times \bb, \Psi)$ to be generated by $h_{\bb}(b_i) = \sigma((n_0, b_i), \dots, (n_{k-1}, b_i))$.  Because these linear orders belong to $\Upsilon[K]$, this is a $K$-embedding.

This a commuting system of maps from the system generating $EM_\tau(I, \Phi)$ to the system generating $EM_\tau(\omega\times I, \Psi)$.  Thus, it generates $h: EM_\tau(I, \Phi) \to EM_\tau( \omega \times I, \Psi)$, which proves the claim.

The proof of the claim above gives a proof of the following result.
\begin{cor}\label{concat-cor}
Suppose that $\Phi, \Psi \in \Upsilon[K]$ and $I$ and $J$ are linear orders so $EM_\tau(I, \Phi) \prec_K EM_\tau(J, \Psi)$.  Then there is $\Psi' \in \Upsilon[K]$ such that, for any linear order $I'$, we have
\begin{itemize}
	\item $EM_\tau(I', \Phi) \prec_K EM_\tau(I', \Psi')$; and
	\item $EM_\tau(I', \Psi) \prec_K EM_\tau(I', \Psi')$.
\end{itemize}
\end{cor}

Third, we show that $M$ has arbitrarily large extensions.  For every linear order $I$, we know that
\begin{itemize}
	\item $\seq{EM_\tau(I, \Psi_i') : i < \delta}$ is $\prec_K$-increasing; and
	\item for every $i < \delta$, $M_i \prec EM_\tau(I, \Psi_i')$.
\end{itemize}

Thus, we can set $M_\delta(I) := \cup_{i < \delta} EM_\tau(I, \Psi_i')$ and we have $M \prec M_\delta(I)$.  Although this takes linear orders to models, it is not clear that this map comes from a blueprint; each $EM(I, \Psi_i')$ could be very different and only relate nicely when the restriction to $\tau$ is taken.  Nonetheless, by taking $I$ large enough, we see that $M$ has an extension of size $\beth_{(2^\lambda)^+}$ and, thus, $\Upsilon[K_{M}]$ is nonempty. \hfill $\dag_{Theorem \ref{main-thm}}$ \\

\begin{cor}
Let $K$ be an AEC.  There is a club of cardinal $\bold{C}$ such that, if $\lambda \in \bold{C}$ and $K_{<\lambda}$ has amalgamation, then every $M \in K_\lambda$ lives inside an $EM$ model.
\end{cor}

{\bf Proof:} Set $\bold{C}:= \{ h_\delta(LS(K)) : \delta$ is limit$\}$.  Then use Theorem \ref{main-thm}.\hfill \dag

Comparing this to Shelah's result in \cite{sh893}, there are some gains and loses.  The main loss is the addition of the amalgamation hypothesis.  Much recent work has taken place under the assumption of amalgamation.  Indeed, one of the most basic tools for analyzing AECs, the Galois type, is robbed of much of its power without amalgamation or some similar replacement.  However, this work also takes place under the assumption of no maximal models, which renders this analysis moot.  On the other hand, inductive analyses might show that amalgamation holds below $\lambda$ and then ask what happens at $K_\lambda$.  

The first main gain is that there is no assumption about the number of models in $K_\lambda$.  As mentioned in the introduction, this means that the conclusion of ``every $M \in K_\lambda$ has $\Upsilon[K_M] \neq \emptyset$'' is unlikely to be a dividing line, at least when $\lambda$ is a sufficiently closed limit.  The second main gain is that there is no cofinality restriction on the cardinality of the size of $M$.

An interesting question is whether the requirement that $\|M\|$ be of a special form (here that it is $h_\delta(LS(K))$ for limit $\delta$) can be removed.  It seems unlikely to the author that maximal models could show up cofinally in an AEC, but a result of this kind does not exist.

\end{document}